\documentclass[12pt,a4paper]{article}

\usepackage{amsmath,amsfonts,amsthm,amssymb}
\usepackage{natbib}
\usepackage{pxfonts}
\usepackage[utf8]{inputenc}
\usepackage[T1]{fontenc}
\usepackage[colorlinks=true,linkcolor=blue, citecolor=blue,urlcolor=blue]{hyperref}
\usepackage{csquotes}
\usepackage{float}
\usepackage{graphicx}
\usepackage{multirow}

\theoremstyle{plain}
\newtheorem{theorem}{Theorem}[section]

\newtheorem{corollary}[theorem]{Corollary}
\newtheorem{proposition}[theorem]{Proposition}

\newtheorem{definition}[theorem]{Definition}

\theoremstyle{definition}

\begin{document}
	
	\title{Asymptotic properties of the normalized discrete associated-kernel estimator for probability mass function}
	\author{Youssef Esstafa\thanks{(\textbf{Corresponding Author}) Le Mans Universit\'e, Laboratoire Manceau de Mathématiques, Avenue Olivier Messiaen, 72085 Le Mans Cedex 09, France.  \href{mailto:Youssef.Esstafa@univ-lemans.fr}{Youssef.Esstafa@univ-lemans.fr}},\; 
		C{\'e}lestin C. Kokonendji\thanks{Universit\'e Bourgogne Franche-Comt{\'e}, Laboratoire de Math{\'e}matiques de Besan\c{c}on UMR 6623 CNRS-UBFC, 16 route de Gray, 25000 Besan{\c c}on, France. \href{mailto:celestin.kokonendji@univ-fcomte.fr}{celestin.kokonendji@univ-fcomte.fr}; \& University of Bangui, Department of Mathematics, B.P. 908 Bangui, Central African Republic. \href{mailto:kokonendji@gmail.com}{kokonendji@gmail.com}}\;  
		and Sobom M. Som\'e\thanks{Universit\'{e} Thomas SANKARA, Laboratoire Sciences et Techniques, 12 BP 417 Ouagadougou 12, Ouagadougou, Burkina Faso.              \href{mailto:sobom.some@uts.bf}{sobom.some@uts.bf}; \& Universit\'{e} Joseph KI-ZERBO, Laboratoire d'Analyse Numérique d'Informatique et de BIOmathématique, 03 BP 7021 Ouagadougou 03,	Ouagadougou, Burkina Faso. \href{mailto:sobom.some@univ-ouaga.bf}{sobom.some@univ-ouaga.bf}}
	}
	\maketitle
	
	\begin{abstract}
		\noindent 
		Discrete kernel smoothing is now gaining importance in nonparametric statistics. In this paper, we investigate some asymptotic properties of the normalized discrete associated-kernel estimator of a probability mass function. We show, under some regularity and non-restrictive assumptions on the associated-kernel, that the normalizing random variable converges in mean square to 1. We then derive the consistency and the asymptotic normality of the proposed estimator. Various families of discrete kernels already  exhibited satisfy the conditions, including the refined CoM-Poisson which is underdispersed and of second-order. Finally, the first-order binomial kernel is discussed and, surprisingly, its normalized  estimator has a suitable asymptotic behaviour through simulations.
	\end{abstract}
	
	
	
	\pagebreak
	
	\section{Introduction}\label{sec:Intro}
	
	The modern notion of a discrete associated kernel for smoothing or estimating discrete functions $f$, defined on the discrete set $\mathbb{T}\subseteq\mathbb{R}$, requires the development of new properties of convergences of the corresponding estimator. The support $\mathbb{T}$ of $f$ is not subject to any restrictive condition; it can be bounded, unbounded, finite or infinite. In this sense, \cite{abdous09} presented some asymptotic properties for non-normalized discrete associated-kernel estimators of a probability mass function (pmf). Several authors pointed out the use of a discrete associated kernel from Dirac and discrete triangular kernels \citep{KSZ07,KZ10} and also from extensions of Dirac kernels  proposed by \cite{Aitchison76} for categorial data and \cite{Wang81}. Furthermore, we have count kernels as the binomial \citep{KSK11} and, recently, the CoM-Poisson \citep{Huang21} kernels which are both underdispersed (\textit{i.e.}, variance less than mean). See also \cite{Harfouche2018} and \cite{Senga17} for other properties. Notice that one can use them to estimate, instead of the pmf, discrete regression or  weighted functions; see, \textit{e.g.}, \cite{KS21}, \cite{Senga14} and, \cite{Senga16}.
	
	Let us firstly fix the refined definition of discrete associated kernel from \cite{KS18} and state in Theorem \ref{ThAbdous09} some important asymptotic properties to be completed in this paper.
	
	\begin{definition}\label{def_DAK} 
		Let $\mathbb{T}\subseteq \mathbb{R}$ be the discrete support of the pmf $f$ to be estimated, $x\in \mathbb{T}$ a target point and $h>0$ a bandwidth. A parameterized pmf $K_{x,h}(\cdot)$ on the discrete support $\mathbb{S}_{x}\subseteq \mathbb{R}$ is called "discrete associated kernel" if the following conditions are satisfied:
		\begin{equation}\label{DAK}
			x\in\mathbb{S}_{x},\;
			\lim_{h\rightarrow 0}\mathbb{E}(Z_{x,h})=x \;and\;
			\lim_{h\rightarrow 0}\mathrm{Var}(Z_{x,h})=\delta\in [0,1),
		\end{equation}
		where $Z_{x,h}$ denotes the discrete random variable with pmf $K_{x,h}(\cdot)$.
	\end{definition}
	
	The choice of the discrete associated kernel referred to as a "second-order" satisfying $\delta=0$ in (\ref{DAK}), ensures the convergence of its corresponding estimator; and, an elementary example is the naive or Dirac kernel for smoothing a very large sample of discrete  data. Otherwise, the convergence of its corresponding estimator is not guaranteed; that is for $\delta=\delta(x)\in (0,1)$ in Definition \ref{def_DAK}, the discrete associated kernel is called of  "first-order" like the well-known binomial kernel. 
	
	Let $X_1,X_2,\ldots, X_n$ be a sample of independent and identically distributed (i.i.d.) discrete random variables having a pmf $f$ on $\mathbb{T}\subseteq\mathbb{R}$.  
	In general, the basical estimator of $f$ is not a pmf. Indeed, for some discrete associated-kernels (\textit{e.g.}, binomial,  triangular and CoM-Poisson), the total mass of the corresponding estimator is not equal to 1. This limit is explained by the fact that the normalizing variable (which is equal to the sum over all the targets belonging to $\mathbb{T}$ of the discrete associated-kernel estimator) was assumed to be equal to 1 only to simplify the calculations. More precisely, one can write both estimators as:
	\begin{equation}\label{eq: normal estim}
		\widehat{f}_n(x)=\frac{\widetilde{f}_n(x)}{C_n}, \;\; x\in \mathbb{T},
	\end{equation}
	with 
	\begin{equation}\label{eq: ftild-Cn}
		\widetilde{f}_n(x)=\frac{1}{n}\sum_{i=1}^nK_{x,h_n}(X_i) \;\;\;\; \text{ and } \;\;\;\; C_n=\sum_{x\in\mathbb{T}}\widetilde{f}_n(x)>0,
	\end{equation}
	where $(h_n)_{n\geq 1}$ is an arbitrary sequence of positive smoothing parameters that satisfies $\lim_{n\to\infty}h_n=0$, while $K_{x,h_n}(\cdot)$ is a suitably chosen discrete kernel function. If $C_n=1$ as for the kernels of Dirac, \cite{Aitchison76} and \cite{Wang81}, one obviously has $\widetilde{f}_n=\widehat{f}_n$. Hence:
	
	\begin{theorem}\citep{abdous09}\label{ThAbdous09} 
		For any $x\in\mathbb{T}$ and under Assumptions (\ref{DAK}) of the second-order (\textit{i.e.}, $\delta=0$), one has 
		\begin{equation*}
			\widetilde{f}_n(x)\xrightarrow[n\to \infty]{L^2,\;a.s.} \ f(x),
		\end{equation*}
		where $` `\overset{L^2,\;a.s.}{\longrightarrow}"$ stands for both ``mean square and almost surely convergences". Furthermore, if $f(x)>0$ then 
		\begin{equation*}
			\left\{\widetilde{f}_n(x)-\mathbb{E}\widetilde{f}_n(x)\right\}\left\{\mathrm{Var}\widetilde{f}_n(x)\right\}^{-1/2}\xrightarrow[n\to \infty]{\mathcal{D}} \mathcal{N}(0,1),
		\end{equation*}
		where $` `\overset{\mathcal{D}}{\longrightarrow}"$ stands for ``convergence in distribution" and $\mathcal{N}(0,1)$ denotes the standard normal distribution.
	\end{theorem}
	
	In this paper we mainly extend Theorem \ref{ThAbdous09} of the non-normalized estimator $\widetilde{f}_n$ of (\ref{eq: ftild-Cn}) to the normalized one $\widehat{f}_n$ of (\ref{eq: normal estim}), introducing new and non-restrictive assumptions with uniformities on the target point in the limits of (\ref{DAK}) and, therefore, changing the types of convergences. As a matter of fact and more importantly, we shall demonstrate the convergence in mean square of the positive normalizing random variable $C_n$ of (\ref{eq: ftild-Cn}) to 1; which clearly completes the similar result in \citet[Theorem 2.1]{KV16}. The following Section \ref{sec:Results} states different  assumptions and their corresponding results with illustrations on the recent CoM-Poisson kernel estimator. The case of the first-order binomial kernel is briefly discussed. Finally, Section \ref{sec:Proofs} is devoted to the detailed proofs.
	
	\section{Results and illustrations}\label{sec:Results}
	
	In order to obtain some soft convergences of the pointwise normalized estimator $\widehat{f}_n(\cdot)$ at $x$, we need quite strong assumptions instead of the most popular  (\ref{DAK}). In this way, we do not use concentration inequalities as in \cite{KV16} as well as \cite{abdous09} through, for instance, an inequality of \cite{Hoeffding63}.
	
	The first set of assumptions is uniformly in the target $x$ and it is satisfied, in our knowledge, by all discrete kernels of Definition \ref{def_DAK} with $\delta=0$:
	$$\mathbf{(A1)}\!:\;\; x\in\mathbb{S}_{x},\;\; \lim_{n\to\infty}\sup_{x\in\mathbb{T}}\left|\mathbb{E}(Z_{x,h_n})-x\right|=0\;\;\mathrm{and}\;\;\lim_{n\to\infty}\sup_{x\in\mathbb{T}}\mathrm{Var}(Z_{x,h_n})=0.$$
	Hence, the following proposition provides a key point to establish the next result on the pointwise probability convergence of $\widehat{f}_n(x)$ defined in (\ref{eq: normal estim}).
	
	\begin{proposition}\label{prop: Cn} Under Assumptions $\mathbf{(A1)}$, the normalizing random variable $C_n$ converges in mean square to 1.
	\end{proposition}
	
	\begin{theorem}[Consistency]\label{ThConsist}
		Under $\mathbf{(A1)}$ and for any $x\in\mathbb{T}$, we have:
		\begin{equation*}
			\widehat{f}_n(x)\xrightarrow[n\to \infty]{\mathbb{P}} f(x),
		\end{equation*}
		where $` `\overset{\mathbb{P}}{\longrightarrow}"$ stands for $` `$convergence in probability$"$. 
	\end{theorem}
	
	On a finite discrete set $\mathbb{T}$, the pointwise and uniform convergences of a sequence of functions are equivalent; hence, the previous results are guaranteed when the discrete associated-kernel satisfies the common set of hypotheses (\ref{DAK}) with $\delta=0$ (see Section \ref{sec:Proofs} for further details).
	\begin{corollary}[Uniform consistency]\label{CorollUC} 
		Suppose that the set $\mathbb{T}$ is finite. Under Conditions (\ref{DAK}) with $\delta=0$, one has
		\begin{equation*}
			\sup_{x\in\mathbb{T}}\left|\widehat{f}_n(x)-f(x)\right|\xrightarrow[n\to \infty]{\mathbb{P}} \ 0.
		\end{equation*}
	\end{corollary}
	
	Regarding to a refined result of the asymptotic normality of $\widehat{f}_n(x)$, it is necessary to quantify the speed of convergence to 0 of $\sup_{x\in\mathbb{T}}\mathrm{Var}(Z_{x,h_n})$ and $\sup_{x\in\mathbb{T}}|\mathbb{E}[Z_{x,h_n}]-x|$ in $\mathbf{(A1)}$. We therefore assume that these two sequences satisfy the  following second set of conditions:
	$$\mathbf{(A2)}\!:\;\; x\in\mathbb{S}_{x},\;\; \sup_{x\in\mathbb{T}}\left|\mathbb{E}(Z_{x,h_n})-x\right|=\mathcal{O}(h_n)\;\;\mathrm{and}\;\;\sup_{x\in\mathbb{T}}\mathrm{Var}(Z_{x,h_n})=\mathcal{O}(h_n).$$
	The previous Assumptions $\mathbf{(A2)}$ and also $\mathbf{(A1)}$ are verified by all the usual kernels of second-order introduced as examples in Section \ref{sec:Intro}. 
	
	\begin{theorem}[Asymptotic normality]\label{th:norm} 
		Let {\bf (A2)} be satisfied.  If the sequence $(h_n)_{n\geq 1}$ is chosen such that $\sqrt{n}h_n\longrightarrow0$ as $n\to\infty$, then, for any $x\in\mathbb{T}$ such that $f(x)>0$, the sequence $\{\sqrt{n}(\widehat{f}_n(x)-f(x))\}_{n\geq 1}$ has a limiting centered normal distribution with variance $f(x)\{1-f(x)\}$.
	\end{theorem}
	
	To conclude this section, we highlight some of our previous results on the recent CoM-Poisson kernel estimator of \cite{Huang21} and compare with the classical binomial one. In fact, we consider the refined version of the CoM-Poisson kernel satisfying \textbf{(A1)} and \textbf{(A2)} as follows: $\mathbb{T}=\mathbb{N}=\mathbb{S}_x$ for each $x\in\mathbb{N}$ and any $h>0$, 
	$$K_{x,h}^{CMP}(z)=\frac{\left\lbrace\lambda(x,1/h)\right\rbrace^z}{(z!)^{1/h}}\left\lbrace D\left(\lambda(x,1/h),1/h\right)\right\rbrace^{-1},
	$$
	where $D(\lambda(x,1/h),1/h)=\sum_{z=0}^\infty[\lambda(x,1/h)]^z/(z!)^{1/h}$ is the normalizing constant and $\lambda:=\lambda(x,1/h)$ represents a function of $x$ and $1/h$ given by the solution of
	\begin{align}\label{condition-COM}
		\sum_{z=0}^\infty\frac{\left\lbrace\lambda(x,1/h)\right\rbrace^z}{(z!)^{1/h}}(z-x)=0.
	\end{align}
This construction implies that $\mathbb{E}(Z_{x,h}^{CMP})=x$ and  
\begin{equation}\label{VarCoMPoiss}
		\mathrm{Var}\left(Z_{x,h}^{CMP}\right)=h\left\lbrace\lambda(x,1/h)\right\rbrace^{h}+\mathcal{O}\left(\left\lbrace\lambda(x,1/h)\right\rbrace^{-h}\right)\;\mathrm{as}\;h\to 0.
	\end{equation}
	
	Indeed, \cite{Huang17} proposed the parametrization via the mean of the original  CoM-Poisson (Conway-Maxwell-Poisson or CMP) distribution; see, \textit{e.g.}, \cite{Shmueli05}, \citet[Section 4.2]{KMB08},  \citet{Gaunt19} and \citet[Section 2.2]{Toledo22} for more details on the original form, asymptotic properties and the relative dispersion with respect to the standard Poisson model. Also demonstrating in the Appendix, the following  proposition points out the mean and the main key of the variance behaviour (\ref{VarCoMPoiss}) of this CoM-Poisson kernel which is of the secnd-order and underdispersed for $h\in (0,1)$.

\begin{proposition}\label{Prop:CMP}
Let $Y$ be a count random variable following  the mean-parametrized CoM-Poisson distribution with location (or mean) parameter $\mu\geq 0$ and dispersion parameter $\nu> 0$ such that its pmf $p(\cdot;\mu,\nu)$ is defined by
	\begin{align}\label{CoM-Poisson}
		p(y;\mu,\nu):=K_{\mu,1/\nu}^{CMP}(y), \ \ \ y\in\mathbb{N}.
	\end{align}
Then $\mathbb{E}(Y)=\mu$ and, when $\{\lambda(\mu,\nu)\}^{1/\nu}\to\infty$ as $\nu\to\infty$, the variance of $Y$ verifies
\begin{equation*}\label{borneCMP}
\mathrm{Var}\left(Y\right)=\frac{1}{\nu}\left[\lambda(\mu,\nu)\right]^{1/\nu}+\mathcal{O}\left(\left\lbrace\lambda(\mu,\nu)\right\rbrace^{-1/\nu}\right).
\end{equation*}
\end{proposition}
	
	As for the binomial kernel of first-order and underdispersed \citep{KSK11}, one has: $\mathbb{T}=\mathbb{N}$, $\mathbb{S}_x=\{0,1,\ldots,x+1\}$ for each $x\in\mathbb{N}$ and $h\in (0,1)$, 
	$$
	K_{x,h}^{B}(z)=\frac{(x+1)!}{z!(x+1-z)!}\left(\frac{x+h}{x+1}\right)^z
	\left(\frac{1-h}{x+1}\right)^{x+1-z}
	$$
	with $\mathbb{E}(Z_{x,h}^{B})=x+h\to x$ as $h\to 0$ and 
	\begin{equation}\label{VarBin}
		\mathrm{Var}\left(Z_{x,h}^{B}\right)=\frac{(x+h)(1-h)}{x+1}.
	\end{equation}
	From Assumptions (\ref{DAK}) and through (\ref{VarBin}), one here has $\delta=\delta(x)=x/(x+1)\in [0,1)$ which does not clearly satisfy the last condition of \textbf{(A1)} as well as for \textbf{(A2)}. Notice that $K_{x,h}^{B}(\cdot)$ is the pmf of the standard binomial distribution $\mathcal{B}(n,p)$ with $n:=x+1$ and $p:=(x+h)/(x+1)$. Nevertheless, we always use the binomial kernel for smoothing count data of small and moderate sample sizes.
	
All numerical studies are here performed using the classical binomial and the recent CoM-Poisson kernel smoothers with the aim to corroborate the previous theoretical results.  Computations have been run on PC 2.30 GHz by using the \textsf{R} software \cite{R20}. Both previous estimators are fitted using the \textsf{Ake} package by \cite{Wansouwe2016} and  the \textsf{mpcmp} one of \cite{Fung20}, respectively. We evaluate the performances of these two discrete associated-kernel estimators with cross-validation choices of the optimal bandwidth parameter. In fact, the optimal bandwidth $h_{cv}$ of $h$  using the cross-validation method is obtained through 
	$$h_{cv} = \arg\min_{h>0}\left[\sum_{x \in \mathbb{T}}{\left\{\widehat{f}_n(x)\right\}^{2}} -\displaystyle\frac{2}{n}\displaystyle\sum_{i=1}^{n}{ \widehat{f}_{n,h,-i}(X_i)}\right],$$ 
	where 
	$$\widehat{f}_{n,h,-i}(X_i)=\displaystyle\frac{1}{n-1}\displaystyle\sum_{\ell=1,\ell\neq i}^n K_{X_{i},h}(X_{\ell})$$
is  being computed as $\widehat{f}_n(X_i)$  without the observation $X_i$. 
		
	We here consider four scenarios which are denoted by A, B, C and D  to simulate count datasets  with respect to the support of both discrete kernels. These scenarios have been considered  to evaluate the performances of both smoothers to deal with zero-inflated, unimodal and multimodal distributions. We shall examine the efficiency of both smoothers via the empirical estimates of $\widehat{C}_n$ and $\widehat{ISE}_n$ of the integrated squared errors (ISE):	
$$\widehat{ISE}_n:=\frac{1}{N_{sim}}\displaystyle\sum_{t=1}^{N_{sim}}\displaystyle\sum_{x \in \mathbb{T}}\left\lbrace\widehat{f}_n(x)- f(x)\right\rbrace^{2} \quad \mbox{ and }\quad  \widehat{C}_n:=\frac{1}{N_{sim}}\displaystyle\sum_{t=1}^{N_{sim}}\displaystyle\sum_{x \in \mathbb{T}}\widetilde{f}_n(x),$$
	where $N_{sim}$ is the number of replications and $n$ corresponds to the sample size which shall be small, medium and large.   
	\begin{itemize}		
		\item Scenario A is generated  by using the Poisson distribution  
		$$f_{A}(x)=\frac{8^{x} e^{-8}}{x!}, \quad x \in \mathbb{N};$$
		\item Scenario B comes from the zero-inflated Poisson distribution 
		$$f_B(x)=\left(\frac{7}{10} \mathbf{1}_{\{x=0\}}\right) + \left(\frac{3}{10}\times\frac{ 10^{x} e^{-10}}{x!}\right), \quad x \in \mathbb{N};$$		
		\item Scenario C is from a mixture of two Poisson distributions 
		$$f_{C}(x)=\left(\frac{2}{5}\times\frac{ 0.5^x e^{-0.5}}{x!}\right)+\left(\frac{3}{5}\times\frac{8^{x} e^{-8}}{x!}\right), \quad x \in \mathbb{N};$$			
		\item  Scenario D comes from a mixture of three Poisson distributions 
		$$f_{D}(x)=\left(\frac{3}{5}\times  \frac{10^x e^{-10} }{x!}\right)+\left(\frac{1}{5}\times \frac{22^{x} e^{-22}}{x!}\right)+\left(\frac{1}{5} \times \frac{50^{x} e^{-50}}{x!}\right), \quad x \in \mathbb{N}.$$		
	\end{itemize}	
	\begin{table}[!htbp]
		\begin{center}
			\caption{Empirical mean values of  $\widehat{C}_n$ and $\widehat{ISE}_n$  with their standard deviations in parentheses over $N_{sim}=100$ replications and with different sample sizes $n=10, 25, 50, 100, 250, 500$ under four Scenarios A, B, C and D by using CoM-Poisson and binomial kernels with cross-validated bandwidth selection.} \label{Tab:C_ISE}\hspace*{-1cm}
			\begin{tabular}{rrrrrrrr}
				\hline \hline
				&$n\;$&\multicolumn{1}{c}{$\widehat{C}_n^{B}$ }&\multicolumn{1}{c}{$\widehat{C}_n^{CMP}$ } &\multicolumn{1}{c}{$\widehat{ISE}_n^{B}$ } &\multicolumn{1}{c}{$\widehat{ISE}_n^{CMP}$ }    \\ \hline 
				\multirow{7}{*}{A}& \multirow{1}{*}{10}   &0.98690 (0.01035)&0.91187 (0.05425) & 0.02961 (0.02706) &0.01466 (0.01863) \\  
			
				&\multirow{1}{*}{25} &0.99321 (0.00452)& 0.96705 (0.02471)& 0.01004 (0.00599) & 0.00861 (0.00999)\\		
			
				&\multirow{1}{*}{50}  &0.99460 (0.00259)&0.98634 (0.01233)& 0.00566 (0.00344) &0.00557 (0.00425) \\
			
				&\multirow{1}{*}{100} & 0.99570 (0.00147) &0.99525 (0.00307) & 0.00271 (0.00204) &0.00291 (0.00268)\\
			
				&\multirow{1}{*}{250} & 0.99685 (0.00076) & 0.99973 (0.00104)& 0.00076 (0.00057)& 0.00131 (0.00105)  \\
				&\multirow{1}{*}{500}&0.99703 (0.00043) &1.00015 (0.00081) &0.00017 (0.00021) & 0.00042 (0.00047) \\
				\hline
			
				\multirow{7}{*}{B}& \multirow{1}{*}{10}   & 0.98663 (0.03655) &0.94781 (0.05083)&0.03326 (0.02456)& 0.02232 (0.01877)\\ 
			
				&\multirow{1}{*}{25} & 0.99788 (0.02049)&0.98053 (0.03101)&0.01392 (0.00871) &0.01054 (0.01033)   \\	
				  
				&\multirow{1}{*}{50}  &1.00869 (0.00826) &1.00095 (0.00739) &0.00696 (0.00330)&0.00573 (0.00350) \\
				&\multirow{1}{*}{100} & 1.01265 (0.00557)& 0.99951 (0.00140)& 0.00352 (0.00204)&0.00343 (0.00227)\\
		
				&\multirow{1}{*}{250} &1.01272 (0.00480) &  0.99921 (0.00111) & 0.00055 (0.00034)&0.00107 (0.00076)\\
				&\multirow{1}{*}{500} 
&1.01460 (0.00238) & 0.99969 (0.00059) &0.00051 (0.00033)& 0.00072 (0.00055) \\
				
				\hline 				
				
				\multirow{7}{*}{C}& \multirow{1}{*}{10}   &0.91176 (0.07860)&1.01341 (0.03870) &0.03842 (0.02421) & 0.02711 (0.02782)\\ 				
				
				&\multirow{1}{*}{25} &0.94838 (0.05030)&1.03786 (0.02800) & 0.01175 (0.00786) &0.01021 (0.00874) \\	
					  
				&\multirow{1}{*}{50}  & 0.98508 (0.02758)& 1.03479 (0.01489) & 0.00499 (0.00297) & 0.00520 (0.00451)\\
				
				&\multirow{1}{*}{100} &1.00242 (0.01048) & 1.02441 (0.01134)&0.00273 (0.00159)&0.00336 (0.00265)\\
				&\multirow{1}{*}{250} & 1.04017 (0.01055) & 1.01232 (0.01244) & 0.00053 (0.00047) & 0.00078 (0.00061) \\
				&\multirow{1}{*}{500}&1.03560 (0.00892) &1.00365 (0.00361)&0.00080 (0.00049)&0.00051 (0.00031)  \\
				\hline 
			
				\multirow{7}{*}{D}& \multirow{1}{*}{10}   &0.99122 (0.00118)&0.95058 (0.01701)&0.02489 (0.00987) &0.01315 (0.02328) \\  
				&\multirow{1}{*}{25} &0.99556 (0.00171)&0.97465 (0.01179) &0.00955 (0.00376) &0.00533 (0.00627)\\		  
				&\multirow{1}{*}{50}  &0.99770 (0.00058)& 0.99276 (0.00479) &0.00296 (0.00146)& 0.00294 (0.00256)\\
				&\multirow{1}{*}{100} & 0.99839 (0.00045)& 0.99711 (0.00209)& 0.00098 (0.00044)&0.00125 (0.00089)\\
				&\multirow{1}{*}{250} 
&0.99889 (0.00018)&0.99919 (0.00080) & 0.00022 (0.00012) & 0.00043 (0.00038)  \\
				&\multirow{1}{*}{500}
& 1.01080 (0.00015) &1.00061 (0.00031) & 0.00020 (0.00010) & 0.00042 (0.00024) \\
				\hline
			\end{tabular}			
		\end{center}
	\end{table} 

Table \ref{Tab:C_ISE} reports some empirical mean values of $\widehat{C}_n$ and $\widehat{ISE}_n$ with their standard deviations using $N_{sim}=100$ replications from Scenarios A, B, C and D to the corresponding sample sizes $n = 10, 25, 50, 100, 250, 500$. For each given subsample and the discrete associated-kernel CoM-Poisson or binomial, we have to compute the related  bandwidth $h_{cv}$ through the cross-validation method before $\widetilde{f}_n$, $\widehat{C}_n$, $\widehat{f}_n$ and finally $\widehat{ISE}_n$. Hence, we observe the following behaviours. Firstly, when the sample size $n$ increases then all standard deviations of Table \ref{Tab:C_ISE} steadly decrease towards $0$. The normalizing constant $\widehat{C}_n^{CMP}$ for CoM-Poisson kernel estimator also becomes more and more precise to $1$ in absolute value; while the  binomial one $\widehat{C}_n^{B}$ moves further away from 1 in absolute value for medium and large sample sizes, in particular for both zero-inflated Scenarios B and C. Next and as expected, the consistent CoM-Poisson smoother is increasingly accurate as sample size increases according to the $\widehat{ISE}_n^{CMP}$ criterion. It is  seemingly better than the binomial one $\widehat{ISE}_n^{B}$, especially for small and moderate sample sizes $n\leq 100$. With enormous surprise and satisfaction, the normalized binomial kernel smoother is also asymptotically consistent in practice, similar to the CoM-Poisson one for all used Scenarios. In fact, this normalizing process of $\widetilde{f}_n$ by $C_n$ for obtaining $\widehat{f}_n$ apparently controls the consistency property of $\widehat{f}_n$, even for the discrete first-order associated-kernel not verifying \textbf{(A1)}. Finally, we can notify the importance of normalization of the discrete associated-kernel estimators of pmfs in practice; see, \textit{e.g.}, \cite{Wansouwe2016} and \cite{KS21} for some illustrations in uni- and multivariate cases.

\begin{figure}[h!]
\hspace{-1.75cm}
\begin{tabular}{ll}
\includegraphics[width=8.9cm, height=8.5cm]{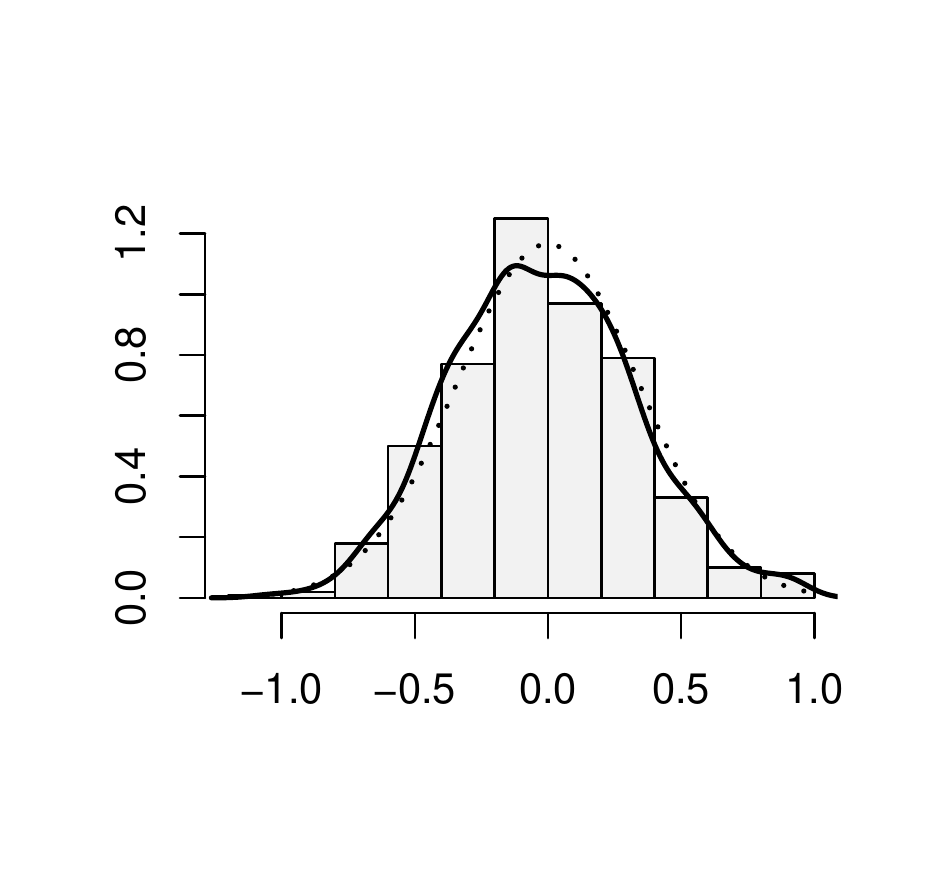}& \hspace{-1.45cm} \includegraphics[width=8.9cm, height=8.5cm]{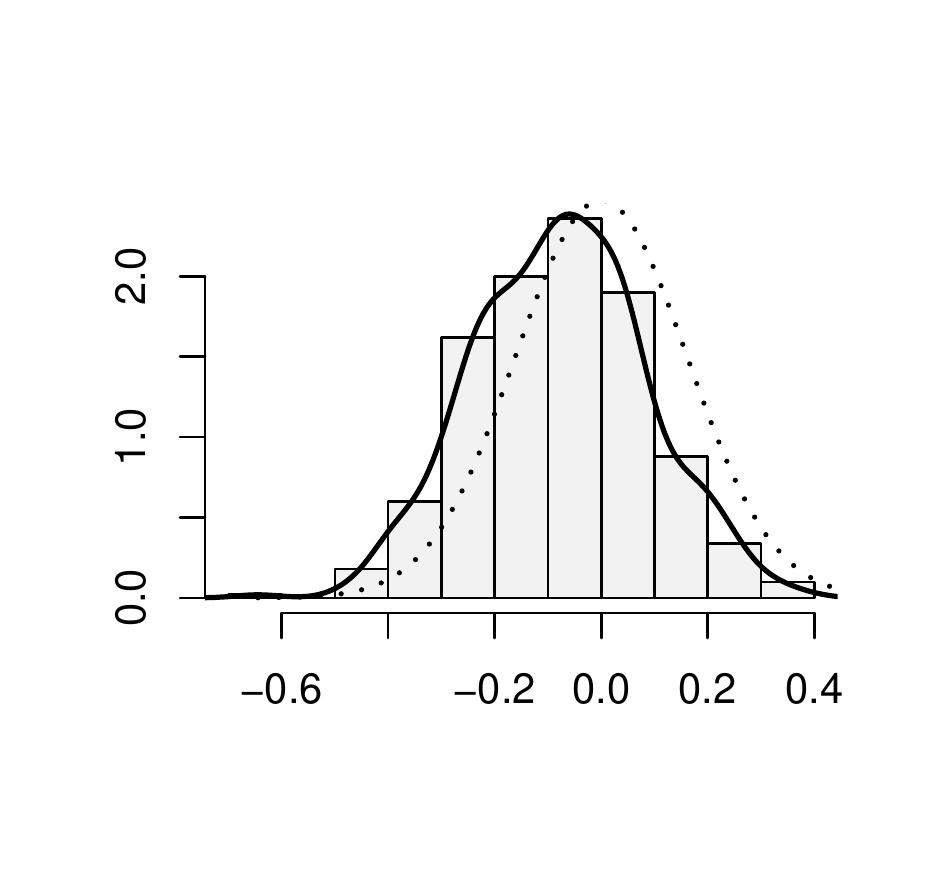}
\end{tabular}
\caption{Empirical distributions of $\sqrt{n}\{\widehat{f}_n(x)-f_A(x)\}$ at $x=6$ over $N_{sim}=500$ independent simulations with $n = 500$ and $h_n=\{\!\sqrt{n}\log(n)\}^{-1}$ using CoM-Poisson (left) and binomial (right) kernels. The smoothed kernel density is displayed in full line, and the centered Gaussian density with the same variance is plotted in dotted line.}\label{Fig:NormAsymp}
\end{figure}

Figure \ref{Fig:NormAsymp} illustrates both empirical distributions of $\sqrt{n}\{\widehat{f}^{CMP}_n(6)-f_A(6)\}$ (left) and $\sqrt{n}\{\widehat{f}^{B}_n(6)-f_A(6)\}$ (right) over $N_{sim}=500$ replications of Scenario A with the sample size $n=500$ and the bandwidth $h_n=\{\!\sqrt{n}\log(n)\}^{-1}$. It is obviously remarkable that the normalized CoM-Poisson kernel estimator is more suitable than the one computed from a binomial kernel, for which the bias increases considerably with the sample size. Once again, these figures confirm the pointwise consistency as well as the pointwise asymptotic normality of the normalized CoM-Poisson kernel estimator, unlike the normalized discrete associated-kernel estimator obtained from the binomial kernel which does not verify our set of hypotheses.
\begin{figure}[!ht]
		\begin{center}
		\resizebox*{14cm}{!}{\includegraphics{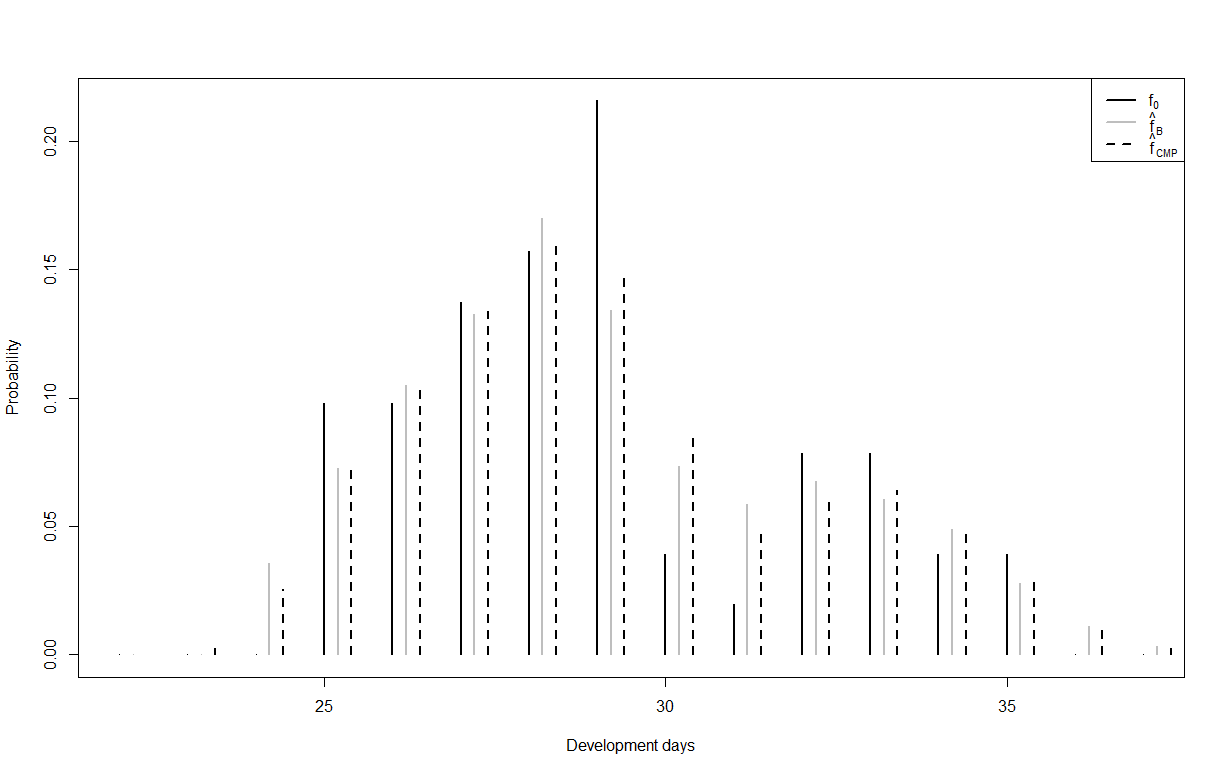}}
			\end{center}
	\caption{Empirical frequency with its corresponding binomial and CoM-Poisson kernel smoothers of count dataset of insect pests on Hura trees with $n=51$; see, \textit{e.g.}, \cite{Senga17}.}
	\label{Fig:hura}
\end{figure}

Concerning an application on real data for pointing out the very competitive CoM-Poisson kernel, both discrete kernel estimators are finally used to smooth a count dataset on development days of insect pests on Hura trees with moderate sample size $n=51$; see \cite{Senga17} and also \cite{Huang21} for applications using these two discrete associated-kernel estimators among others. Practical performances are here examined via the cross-validation method and the empirical criterion of $ISE$: $\widehat{ISE}_0:=\sum_{x\in\mathbb{T}\subseteq\mathbb{N}}[\widehat{f}_{n}(\mathbf{x})-f_0(\mathbf{x})]^2$, where  $f_0(\cdot)$ is the empirical or naive estimator. The CoM-Poisson kernel appears to be the best with $h_{cv}^{CMP}=0.01865$, $\widehat{C}_n^{CMP}=0.99994$ and $\widehat{ISE}_0^{CMP}=0.00967$ followed by the binomial smoother with $h_{cv}^{B}=0.00601$,  $\widehat{C}_n^B=0.99913$ and finally $\widehat{ISE}_0^B=0.01232$; see Figure \ref{Fig:hura} for graphical representations. Notice that, for the same dataset, \cite{Senga17} produced $h_{cv}^{B}=0.02$ and $\widetilde{ISE}_0^B=0.0104$ for the non-normalized binomial estimation $\widetilde{f}_{n}^B$. While \cite{Huang21} only presented $h_{cv}^{CMP}=0.0251$ without $\widetilde{ISE}_0^{CMP}$ for the non-normalized CoM-Poisson estimation $\widetilde{f}_{n}^{CMP}$ with precisions to fit a non-zero probability outside of the observed range and also to preserve the sample mean of the dataset.

\section{Proofs of results}\label{sec:Proofs}
	
	\begin{proof}[Proof of Proposition \ref{prop: Cn}] 
		Firstly, one easily has the following decomposition: 
		\begin{equation}\label{eq: L2Cn}
			\mathbb{E}\left[\left|C_n-1\right|^2\right]=\mathrm{Var}(C_n)+\left(\mathbb{E}[C_n]-1\right)^2.
		\end{equation}
		We use Equation \eqref{eq: ftild-Cn} and the fact that the $X_i$'s are i.i.d. to obtain
		\begin{align}\label{eq: VarCn}
			\mathrm{Var}(C_n)&=\sum_{x\in\mathbb{T}}\sum_{y\in\mathbb{T}}\mathrm{Cov}\left(\widetilde{f}_n(x),\widetilde{f}_n(y)\right)\nonumber\\
			&=\frac{1}{n^2}\sum_{x\in\mathbb{T}}\sum_{y\in\mathbb{T}}\sum_{i=1}^n\sum_{j=1}^n\mathrm{Cov}\left(K_{x,h_n}(X_i),K_{y,h_n}(X_j)\right)\nonumber\\
			&=\frac{1}{n^2}\sum_{x\in\mathbb{T}}\sum_{y\in\mathbb{T}}\sum_{i=1}^n\mathrm{Cov}\left(K_{x,h_n}(X_i),K_{y,h_n}(X_i)\right)\nonumber\\
			&=\frac{1}{n}\sum_{x\in\mathbb{T}}\sum_{y\in\mathbb{T}}\mathrm{Cov}\left(K_{x,h_n}(X_1),K_{y,h_n}(X_1)\right)\nonumber\\
			&=\frac{1}{n}\sum_{x\in\mathbb{T}}\mathrm{Var}\left[K_{x,h_n}(X_1)\right]+\frac{1}{n}\sum_{x\in\mathbb{T}}\sum_{y\in\mathbb{T}\setminus\{x\}}\mathrm{Cov}\left(K_{x,h_n}(X_1),K_{y,h_n}(X_1)\right).
		\end{align}
		The bias term in \eqref{eq: L2Cn} can be explicitly rewritten as:
		\begin{align}\label{eq:bias}
			\mathbb{E}[C_n]-1&=\sum_{x\in\mathbb{T}}\mathbb{E}[\widetilde{f}_n(x)]-\sum_{x\in\mathbb{T}}f(x)\nonumber\\
			&=\sum_{x\in\mathbb{T}}\sum_{z\in\mathbb{T}\cap\mathbb{S}_x}K_{x,h_n}(z)f(z)-\sum_{x\in\mathbb{T}}\sum_{z\in\mathbb{S}_x}K_{x,h_n}(z)f(x)\nonumber\\
			&=E_{1,n}-E_{2,n},
		\end{align}
		with 
		\begin{align*}
			E_{1,n}&=\sum_{x\in\mathbb{T}}\sum_{z\in\mathbb{T}\cap\mathbb{S}_x}\left(f(z)-f(x)\right)K_{x,h_n}(z)
		\end{align*}
		and
		\begin{align*}
			E_{2,n}&=\sum_{x\in\mathbb{T}}f(x)\sum_{z\in\overline{\mathbb{T}}\cap\mathbb{S}_x}K_{x,h_n}(z),
		\end{align*}
		where $\overline{\mathbb{T}}$ is the set $\mathbb{R}\setminus\mathbb{T}$.
		
		Secondly, to make the proof more readable, we divide it in two steps according to the covergences to 0 of both variance and bias terms in \eqref{eq: L2Cn}.
				
		\paragraph{$\diamond$ Step 1: Convergence to 0 of the variance term in \eqref{eq: L2Cn}.}
		Under $\mathbf{(A1)}$, one can prove that the first term in the right hand side of \eqref{eq: VarCn} converges to 0. As a matter of fact, observe first that
		\begin{align}\label{eq: VarK}
			\mathrm{Var}\left[K_{x,h_n}(X_1)\right]-f(x)\left\lbrace 1-f(x)\right\rbrace &=\left\lbrace\sum_{z\in\mathbb{T}\cap\mathbb{S}_x}(K_{x,h_n}(z))^2f(z)-f(x)\right\rbrace\nonumber\\
			&\quad +\left\lbrace f(x)^2-\left(\mathbb{E}\left[\widetilde{f}_n(x)\right]\right)^2\right\rbrace\nonumber\\
			&=:F_{1,n}(x)+F_{2,n}(x).
		\end{align}
		
		The sets $\mathbb{T}$ and $\mathbb{S}_x$ are discrete. So, one can find a finite constant $\alpha>0$ (which does not depend on $x$) such that $|z-x|\geq \alpha$ for any $z$ in $\mathbb{T}\cup\mathbb{S}_x\setminus\{x\}$. Hence, the use of Markov's inequality and Assumptions $\mathbf{(A1)}$ lead to deduce that the first sequence  $(F_{1,n})_{n\geq 1}$ in (\ref{eq: VarK}) converges uniformly on $\mathbb{T}$ to 0 as follows:
		\begin{align*}
			\sup_{x\in\mathbb{T}}\left|F_{1,n}(x)\right|&\leq \sup_{x\in\mathbb{T}} \Bigg\{ f(x)\left(1-K_{x,h_n}(x)\right)K_{x,h_n}(x)\\
			&\left.\quad+\sum_{z\in\mathbb{T}\cap\mathbb{S}_x\setminus \{x\}}\left|f(z)K_{x,h_n}(z)-f(x)\right|K_{x,h_n}(z)+f(x)\sum_{z\in\overline{\mathbb{T}}\cap\mathbb{S}_x}K_{x,h_n}(z)\right\rbrace\\
			&\leq 4 \sup_{x\in\mathbb{T}}\left\lbrace 1-K_{x,h_n}(x)\right\rbrace\\
			&\leq 4 \sup_{x\in\mathbb{T}}\mathbb{P}\left(\left|Z_{x,h_n}-x\right|\geq \alpha\right)\\
			&\leq \frac{4}{\alpha^2} \left\lbrace\sup_{x\in\mathbb{T}}\mathrm{Var}\left(Z_{x,h_n}\right)+\left(\sup_{x\in\mathbb{T}}\left|\mathbb{E}\left[Z_{x,h_n}\right]-x\right|\right)^2\right\rbrace\to 0, \text{ as } n\to\infty.
		\end{align*}
		
		Similarly, we show that the second sequence $(F_{2,n})_{n\geq 1}$ in (\ref{eq: VarK}) also  converges uniformly on $\mathbb{T}$ to 0. Let us be more precise. Note that
		\begin{align*}
			\left|F_{2,n}(x)\right|&=f(x)\left|\mathbb{E}\left[\widetilde{f}_n(x)\right]-f(x)\right|+\mathbb{E}\left[\widetilde{f}_n(x)\right]\left|\mathbb{E}\left[\widetilde{f}_n(x)\right]-f(x)\right|\\
			&\leq 2 \left|\mathbb{E}\left[\widetilde{f}_n(x)\right]-f(x)\right|\\
			&\leq 2 \left\lbrace\sum_{z\in\mathbb{T}\cap\mathbb{S}_x\setminus{\{x\}}}K_{x,h_n}(z)\left|f(z)-f(x)\right|+f(x)\sum_{z\in\overline{\mathbb{T}}\cap\mathbb{S}_x}K_{x,h_n}(z)\right\rbrace\\
			&\leq 6 \mathbb{P}\left(Z_{x,h_n}\ne x\right),
		\end{align*}
		where we have used the fact that if $z\in\overline{\mathbb{T}}\cap\mathbb{S}_x$, then necessarily $z\ne x$.
		Arguing as before, we obtain the expected uniform convergence of $(F_{2,n})_{n\geq 1}$ to 0. 
		
		Consequently, from Equation \eqref{eq: VarK} 
		the sequence $\{\mathrm{Var}(K_{\cdot,h_n}(X_1))\}_{n\geq 1}$ converges uniformly on $\mathbb{T}$ to $\ell_1: x\in\mathbb{T}\mapsto f(x)\{1-f(x)\}\in\mathbb{R}$. Finally, there exists $N\in\mathbb{N}^*$ such that for any $n\geq N$, we have 
		\begin{equation*}
			\frac{1}{n}\sum_{x\in\mathbb{T}}\mathrm{Var}\left[K_{x,h_n}(X_1)\right]\leq \frac{2}{n}\sum_{x\in\mathbb{T}}f(x)\{1-f(x)\}\leq \frac{2}{n},
		\end{equation*}
		and the enacted convergence is obtained.
		
		We now deal with the second term on the right hand side of Equation \eqref{eq: VarCn}. 
		Using the definition of the non-normalized associated-kernel estimator $\widetilde{f}_n(\cdot)$ introduced in Equation \eqref{eq: ftild-Cn}, one can write that for any $x,y\in\mathbb{T}$ and all $n\geq 1$,
		\begin{align*}
			\mathrm{Cov}\left(K_{x,h_n}(X_1),K_{y,h_n}(X_1)\right)&=\sum_{z\in\mathbb{T}\cap\mathbb{S}_x\cap\mathbb{S}_y}K_{x,h_n}(z)K_{y,h_n}(z)f(z)\\
			&\quad-\left(\sum_{z\in\mathbb{T}\cap\mathbb{S}_x}K_{x,h_n}(z)f(z)\right)\left(\sum_{z\in\mathbb{T}\cap\mathbb{S}_y}K_{y,h_n}(z)f(z)\right)\\
			&=\sum_{z\in\mathbb{T}\cap\mathbb{S}_x\cap\mathbb{S}_y}K_{x,h_n}(z)K_{y,h_n}(z)f(z)-\mathbb{E}\left[\widetilde{f}_n(x)\right]\mathbb{E}\left[\widetilde{f}_n(y)\right].
		\end{align*}
		It then follows that 
		\begin{align*}
			\mathrm{Cov}\left(K_{x,h_n}(X_1),K_{y,h_n}(X_1)\right)+f(x)f(y)&=\sum_{z\in\mathbb{T}\cap\mathbb{S}_x\cap\mathbb{S}_y}K_{x,h_n}(z)K_{y,h_n}(z)f(z)\\
			&\quad-\left(\mathbb{E}\left[\widetilde{f}_n(x)\right]-f(x)\right)\mathbb{E}\left[\widetilde{f}_n(y)\right]\\
			&\qquad-\left(\mathbb{E}\left[\widetilde{f}_n(y)\right]-f(y)\right)f(x).
		\end{align*}
		Thus, one has 
		\begin{equation*}
			\sup_{(x,y)\in\mathbb{T}^2\atop x\ne y}\left|\mathrm{Cov}\left(K_{x,h_n}(X_1),K_{y,h_n}(X_1)\right)+f(x)f(y)\right|\leq G_{1,n}+2G_{2,n},
		\end{equation*}
		with 
		\begin{equation*}
			G_{1,n}=\sup_{(x,y)\in\mathbb{T}^2\atop x\ne y}\sum_{z\in\mathbb{T}\cap\mathbb{S}_x\cap\mathbb{S}_y}K_{x,h_n}(z)K_{y,h_n}(z)f(z) \ \ \ \text{ and } \ \ \ G_{2,n}=\sup_{x\in\mathbb{T}}\left|\mathbb{E}\left[\widetilde{f}_n(x)\right]-f(x)\right|.
		\end{equation*}
		Following the same arguments used to prove the convergence of $F_{2,n}(x)$ to 0, we show that $(G_{2,n})_{n\geq 1}$ converges to 0. Indeed, observe that
		\begin{align*}
			\sum_{z\in\mathbb{T}\cap\mathbb{S}_x\cap\mathbb{S}_y}K_{x,h_n}(z)K_{y,h_n}(z)f(z)&=K_{x,h_n}(x)K_{y,h_n}(x)f(x)\\
			&\quad +\sum_{z\in\mathbb{T}\cap\mathbb{S}_x\cap\mathbb{S}_y\atop z\ne x}K_{x,h_n}(z)K_{y,h_n}(z)f(z)\\
			&\leq \sum_{z\in\mathbb{S}_y\setminus\{y\}}K_{y,h_n}(z)+\sum_{z\in\mathbb{S}_x\setminus\{x\}}K_{x,h_n}(z)\\
			&=\mathbb{P}\left(Z_{y,h_n}\ne y\right)+\mathbb{P}\left(Z_{x,h_n}\ne x\right).
		\end{align*}
		Hence, the same lines of proof given before can be reproduced to show that $(G_{1,n})_{n\geq 1}$ converges to 0. The sequence $\{\mathrm{Cov}\left(K_{\cdot,h_n}(X_1),K_{\cdot,h_n}(X_1)\right)\}_{n\geq 1}$ is uniformly convergent on the set $\Delta=\{(x,y)\in\mathbb{T}^2, x\ne y\}$ with limit $\ell_2:(x,y)\in \Delta\mapsto -f(x)f(y)\in\mathbb{R}$. It can therefore be easily demonstrated that there exists a positive integer $\widetilde{N}$ such that for any $n\geq\widetilde{N}$, one has 
		\begin{align*}
			\frac{1}{n}\sum_{x\in\mathbb{T}}\sum_{y\in\mathbb{T}\setminus\{x\}}\mathrm{Cov}\left(K_{x,h_n}(X_1),K_{y,h_n}(X_1)\right)&\leq \frac{1}{n}\sum_{x\in\mathbb{T}}f(x)\sum_{y\in\mathbb{T}\setminus\{x\}}f(y)\leq\frac{1}{n}.
		\end{align*}
		This completes the proof of the convergence to 0 of the second term on the right hand side of \eqref{eq: VarCn} and, finally, the proof of the convergence to 0 of the variance term in \eqref{eq: L2Cn}.
		
		\paragraph{$\diamond$ Step 2: Convergence to 0 of the bias term in \eqref{eq: L2Cn}.} 
		The sequence $(E_{2,n})_{n\geq 1}$ introduced in Equation \eqref{eq:bias} clearly converges to 0 since 
		\begin{equation*}
			E_{2,n}\leq \sup_{x\in\mathbb{T}}\mathbb{P}\left(Z_{x,h_n}\ne x\right)\to 0, \text{ as } n\to\infty.
		\end{equation*}
		
		We now use Equation \eqref{eq: ftild-Cn} and the same arguments developed in the previous step to obtain that
		\begin{align*}
			\left|E_{1,n}\right|&\leq\sum_{x\in\mathbb{T}}\left|\left(\sum_{z\in\mathbb{T}\cap\mathbb{S}_x}f(z)K_{x,h_n}(z)\right)-f(x)+f(x)\left\lbrace 1-K_{x,h_n}(x)\right\rbrace\right|\\
			&\quad+\sum_{x\in\mathbb{T}}f(x)\sum_{z\in\mathbb{T}\cap\mathbb{S}_x\setminus\{x\}}f(z)K_{x,h_n}(z)\\
			&\leq \sum_{x\in\mathbb{T}}\left|\mathbb{E}\left[\widetilde{f}_n(x)-f(x)\right]\right|+2\sup_{x\in\mathbb{T}}\mathbb{P}\left(Z_{x,h_n}\ne x\right)\\
			&\leq \sup_{x\in\mathbb{T}}\frac{\sup_{x\in\mathbb{T}}\left|\mathbb{E}\left[\widetilde{f}_n(x)-f(x)\right]\right|}{f(x)}+2\sup_{x\in\mathbb{T}}\mathbb{P}\left(Z_{x,h_n}\ne x\right)\to 0, \text{ as } n\to\infty.
		\end{align*}
		
		Consequently, the bias term in \eqref{eq: L2Cn} converges to 0. This concludes the proof of the proposition.
	\end{proof}
	
	\begin{proof}[Proof of Theorem \ref{ThConsist}] 
		Note that, for any $x\in\mathbb{T}$, one can express
		\begin{align*}
			\widehat{f}_n(x)-f(x)&=\frac{1}{C_n}\left\lbrace\left(\widetilde{f}_n(x)-f(x)\right)+(1-C_n)f(x)\right\rbrace.
		\end{align*}
		Theorem \ref{ThAbdous09} of \cite{abdous09} recalls that $\widetilde{f}_n(x)$ converges in mean square to $f(x)$; such result obviously remains valid in our context. Proposition \ref{prop: Cn} and Slutsky's theorem complete the proof.
	\end{proof}
	
	\begin{proof}[Proof of Corollary \ref{CorollUC}]
		It is enough to observe that 
		\begin{align*}
			\sup_{x\in\mathbb{T}}\left|\widehat{f}_n(x)-f(x)\right|&=\frac{1}{C_n}\sup_{x\in\mathbb{T}}\left|\widetilde{f}_n(x)-f(x)+f(x)\left(1-C_n\right)\right|\\
			&\leq \frac{1}{C_n}\left\lbrace\left|1-C_n\right|+\sup_{x\in\mathbb{T}}\left|\widetilde{f}_n(x)-f(x)\right|\right\rbrace.
		\end{align*}
		Consequently, Proposition \ref{prop: Cn}, Theorem \ref{ThConsist} and the continuous mapping theorem easily allow to deduce the corollary.
	\end{proof}
	
	\begin{proof}[Proof of Theorem \ref{th:norm}] 
		From the end of Theorem \ref{ThAbdous09}, one may first observe that
		\begin{align}\label{eq:norm}
			\sqrt{n}\left(\widehat{f}_n(x)-f(x)\right)&=\frac{1}{C_n}\left(\frac{\widetilde{f}_n(x)-\mathbb{E}\left[\widetilde{f}_n(x)\right]}{\sqrt{\mathrm{Var}\left\lbrace\widetilde{f}_n(x)\right\rbrace}}\right)\sqrt{n\mathrm{Var}\left\lbrace\widetilde{f}_n(x)\right\rbrace}\nonumber\\
			&\quad+\frac{\sqrt{n}}{C_n}\left(\mathbb{E}\left[\widetilde{f}_n(x)\right]-C_nf(x)\right)\nonumber\\
			&=\frac{1}{C_n}\left(\frac{\widetilde{f}_n(x)-\mathbb{E}\left[\widetilde{f}_n(x)\right]}{\sqrt{\mathrm{Var}\left\lbrace\widetilde{f}_n(x)\right\rbrace}}\right)\sqrt{n\mathrm{Var}\left\lbrace\widetilde{f}_n(x)\right\rbrace}\nonumber\\
			&\quad+\frac{\sqrt{n}}{C_n}\left(\mathbb{E}\left[\widetilde{f}_n(x)\right]-f(x)\right)+\frac{\sqrt{n}}{C_n}(1-C_n)f(x).
		\end{align}
		
		Let $(Y_{n,i})$ be the rowwise i.i.d. triangular array defined by
		\begin{equation*}
			Y_{n,i}=\frac{K_{x,h_n}(X_i)-\mathbb{E}\left[K_{x,h_n}(X_i)\right]}{\sqrt{n\mathrm{Var}\left( K_{x,h_n}(X_i)\right)}}, \ \ \ i=1,\dots,n.
		\end{equation*}
		It is clear that for all $n\geq 1$ and any $i=1,\dots,n$, 
		$$\mathbb{E}\left[Y_{n,i}\right]=0 \ \ \ \text{ and } \ \ \ \sum_{i=1}^n\mathbb{E}\left[Y_{n,i}^2\right]=1.$$
		Moreover, since $\mathrm{Var}[K_{x,h_n}(X_1)]\longrightarrow f(x)(1-f(x))$ as $n\to\infty$ (see Step 1 of the proof of Proposition \ref{prop: Cn}) and using the fact that the $X_i$'s are i.i.d.,  one has
		\begin{align*}
			\sum_{i=1}^n\mathbb{E}\left[\left|Y_{n,i}\right|^3\right]&=\frac{\mathbb{E}\left[\left|K_{x,h_n}(X_1)-\mathbb{E}\left[K_{x,h_n}(X_1)\right]\right|^3\right]}{\sqrt{n}\left\lbrace\mathrm{Var}\left[ K_{x,h_n}(X_1)\right]\right\rbrace^{3/2}}\\
			&\leq \frac{1}{\sqrt{n}\left\lbrace\mathrm{Var}\left[ K_{x,h_n}(X_1)\right]\right\rbrace^{3/2}}\longrightarrow 0, \text{ as } n\to\infty.
		\end{align*}
		Thus, Lindeberg's Theorem implies that 
		\begin{equation*}
			\frac{\widetilde{f}_n(x)-\mathbb{E}\left[\widetilde{f}_n(x)\right]}{\sqrt{\mathrm{Var}\left\lbrace\widetilde{f}_n(x)\right\rbrace}}=\sum_{i=1}^nY_{n,i}\xrightarrow[n\to \infty]{\mathbb{\mathcal{D}}} \ \mathcal{N}(0,1).
		\end{equation*}
Note also that 
		\begin{equation*}
			\sqrt{n\mathrm{Var}\left\lbrace\widetilde{f}_n(x)\right\rbrace}=\sqrt{\mathrm{Var}\left[ K_{x,h_n}(X_1)\right]}\longrightarrow \sqrt{f(x)\left\lbrace 1-f(x)\right\rbrace}, \ \ \text{ as } \ \ n\to\infty.
		\end{equation*}
		In view of Proposition \ref{prop: Cn}, Equation \eqref{eq:norm} and Slutsky's theorem, it is therefore sufficient to prove that the two sequences $\{\sqrt{n}(\mathbb{E}[\widetilde{f}_n(x)]-f(x))\}_{n\geq 1}$ and $\{\sqrt{n}(1-C_n)\}_{n\geq 1}$ converge in probability to 0 for obtaining the expected convergence in distribution stated in Theorem \ref{th:norm}.
		
		Indeed, we consider similar arguments as in the proof of the uniform convergence of $F_{2,n}$ on $\mathbb{T}$ (see Step 1 of the proof of Proposition \ref{prop: Cn}) and the assumptions of Theorem \ref{th:norm} to deduce that 
		\begin{equation*}
			\sqrt{n}\left|\mathbb{E}\left[\widetilde{f}_n(x)\right]-f(x)\right|\leq 3\sqrt{n}\sup_{x\in\mathbb{T}}\mathbb{P}\left(Z_{x,h_n}\ne x\right)\longrightarrow 0, \text{ as } n\to\infty.
		\end{equation*}
		To complete the proof of the theorem, note that by using the results of Step 1 in the proof of Proposition \ref{prop: Cn}, there exists $N\in\mathbb{N}^*$ such that for any $n\geq N$, we have 
		\begin{align*}
			\mathrm{Var}\left[\sqrt{n}(1-C_n)\right]&=\sum_{x\in\mathbb{T}}\left\lbrace\mathrm{Var}\left[K_{x,h_n}(X_1)\right]-f(x)\left\lbrace 1-f(x)\right\rbrace\right\rbrace\\
			&\quad +\sum_{x\in\mathbb{T}}\sum_{y\in\mathbb{T}\setminus \{x\}}\left\lbrace\mathrm{Cov}\left(K_{x,h_n}(X_1),K_{y,h_n}(X_1)\right)+f(x)f(y)\right\rbrace\\
			&\leq \kappa^*h_n\longrightarrow 0, \text{ as } n\to\infty,
		\end{align*}
		where $\kappa^*$ is a positive constant.
		
		Similarly, one can use the results of Step 2 in the proof of Proposition \ref{prop: Cn} to prove
		\begin{equation*}
			\mathbb{E}\left[\sqrt{n}(1-C_n)\right]\leq \kappa^{**}\sqrt{n}h_n\longrightarrow 0, \text{ as } n\to\infty,
		\end{equation*}
		where $\kappa^{**}$ is a positive constant. This completes the proof of the theorem.
	\end{proof}
	
	\begin{proof}[Proof of Proposition \ref{Prop:CMP}] 
From Equations \eqref{condition-COM} and \eqref{CoM-Poisson}, one firstly has
\begin{align*}
\mathbb{E}\left[Y\right]&=\mu+\frac{1}{D\left(\lambda(\mu,\nu),\nu\right)}\sum_{y=0}^\infty(y-\mu)\frac{\left\lbrace\lambda(\mu,\nu)\right\rbrace^y}{(y!)^{\nu}}=\mu.
\end{align*}
Secondly, one can observe that
\begin{align}\label{varCMP1}
\lambda(\mu,\nu)\frac{\partial\mathbb{E}[Y]}{\partial\lambda(\mu,\nu)}&=\frac{1}{D\left(\lambda(\mu,\nu),\nu\right)} \sum_{y=0}^\infty y^2\frac{\left\lbrace\lambda(\mu,\nu)\right\rbrace^y}{(y!)^{\nu}}\nonumber\\
&\quad+\lambda(\mu,\nu)\frac{\partial \left\lbrace 1/ D\left(\lambda(\mu,\nu),\nu\right) \right\rbrace}{\partial\lambda(\mu,\nu)}\sum_{y=0}^\infty y\frac{\left\lbrace\lambda(\mu,\nu)\right\rbrace^y}{(y!)^{\nu}}\nonumber\\
&=\mathbb{E}[Y^2]+\lambda(\mu,\nu)D\left(\lambda(\mu,\nu),\nu\right)\frac{\partial \left\lbrace 1/ D\left(\lambda(\mu,\nu),\nu\right) \right\rbrace}{\partial\lambda(\mu,\nu)}\mathbb{E}\left[Y\right].
\end{align}
Moreover, a direct calculation of derivative leads us to
\begin{align}\label{varCMP2}
\lambda(\mu,\nu)D\left(\lambda(\mu,\nu),\nu\right)\frac{\partial \left\lbrace 1/ D\left(\lambda(\mu,\nu),\nu\right) \right\rbrace}{\partial\lambda(\mu,\nu)}&=-\lambda(\mu,\nu)\frac{\partial D\left(\lambda(\mu,\nu),\nu\right) / \partial\lambda(\mu,\nu)}{D\left(\lambda(\mu,\nu),\nu\right) }\nonumber\\
&=\frac{-1}{D\left(\lambda(\mu,\nu),\nu\right)}\sum_{y=0}^\infty y\frac{\left\lbrace\lambda(\mu,\nu)\right\rbrace^y}{(y!)^{\nu}}\nonumber\\
&=-\mathbb{E}\left[Y\right].
\end{align}
Hence, using \eqref{varCMP1} and \eqref{varCMP2} we deduce 
\begin{equation*}
\mathrm{Var}\left(Y\right)=\lambda(\mu,\nu)\frac{\partial\mathbb{E}[Y]}{\partial\lambda(\mu,\nu)}.
\end{equation*}
We still consider Equation \eqref{condition-COM} to obtain
\begin{align*}
\lambda(\mu,\nu)=\mu\left(\frac{D\left(\lambda(\mu,\nu),\nu\right)}{\partial D\left(\lambda(\mu,\nu),\nu\right)/ \partial\lambda(\mu,\nu)}\right),
\end{align*} 
which implies that
\begin{align}\label{vardlogD}
\mathrm{Var}\left(Y\right)&=\lambda(\mu,\nu)\frac{\partial}{\partial\lambda(\mu,\nu)}\left[\lambda(\mu,\nu)\frac{\partial\log D\left(\lambda(\mu,\nu),\nu\right)}{\partial\lambda(\mu,\nu)}\right]\nonumber\\
&=\lambda(\mu,\nu)\frac{\partial\log D\left(\lambda(\mu,\nu),\nu\right)}{\partial\lambda(\mu,\nu)}+\left\lbrace\lambda(\mu,\nu)\right\rbrace^2\frac{\partial^2\log D\left(\lambda(\mu,\nu),\nu\right)}{\left\lbrace\partial\lambda(\mu,\nu)\right\rbrace^2}.
\end{align}
According to \citet[Theorem 1]{Gaunt19} under $\{\lambda(\mu,\nu)\}^{1/\nu}\to\infty$ as $\nu\to\infty$, one can write
\begin{align*}
\log D\left(\lambda(\mu,\nu),\nu\right)&=\nu\left\lbrace\lambda(\mu,\nu)\right\rbrace^{1/\nu} - \frac{\nu-1}{2\nu}\log\lambda(\mu,\nu) -\log\left\lbrace (2\pi)^{(\nu-1)/2}\sqrt{\nu}\right\rbrace\\
&\quad+\frac{\nu^2-1}{24\nu}\left\lbrace\lambda(\mu,\nu)\right\rbrace^{-1/\nu}+\frac{\nu^2-1}{48\nu^2}\left\lbrace\lambda(\mu,\nu)\right\rbrace^{-2/\nu}+\mathcal{O}\left(\left\lbrace\lambda(\mu,\nu)\right\rbrace^{-3/\nu}\right);
\end{align*}
and, therefore, one deduces the following two terms of \eqref{vardlogD} as:
$$
\lambda(\mu,\nu)\frac{\partial\log D\left(\lambda(\mu,\nu),\nu\right)}{\partial\lambda(\mu,\nu)}=
\left\lbrace\lambda(\mu,\nu)\right\rbrace^{1/\nu}-\frac{\nu-1}{2\nu}+\mathcal{O}\left(\left\lbrace\lambda(\mu,\nu)\right\rbrace^{-1/\nu}\right)
$$
and 
$$\left\lbrace\lambda(\mu,\nu)\right\rbrace^2\frac{\partial^2\log D\left(\lambda(\mu,\nu),\nu\right)}{\left\lbrace\partial\lambda(\mu,\nu)\right\rbrace^2}= \frac{1-\nu}{\nu}\left\lbrace\lambda(\mu,\nu)\right\rbrace^{1/\nu}+\frac{\nu-1}{2\nu}+\mathcal{O}\left(\left\lbrace\lambda(\mu,\nu)\right\rbrace^{-1/\nu}\right).
$$
This is enough to complete the proof.
\end{proof}

\end{document}